\documentclass[12pt,reqno]{amsart}
\usepackage{amsfonts,amsmath,amssymb}

\allowdisplaybreaks
\advance\textheight by \topskip

\usepackage{mathrsfs}

\usepackage[latin1]{inputenc}
\usepackage{graphicx}

\usepackage{tikz}
\newtheorem{theorem}{Theorem}
\usepackage{euscript}

\def\[{[\! [}
\def\]{]\! ]}

\begin{document}

\title{An online bin-packing problem with an underlying ternary structure}

\author[H.~Prodinger]{Helmut Prodinger}

\address{Helmut Prodinger,
Mathematics Department, Stellenbosch University,
7602 Stellenbosch, South Africa.}
\email{hproding@sun.ac.za}

\date{\today}


\begin{abstract}
	Following an orginal idea by Kn\"odel, an online bin-packing problem is considered where the the large
	items arrive in double-packs. The dual problem where the small items arrive in double-packs is
	also considered. The enumerations have a ternary random walk flavour, and for the enumeration,
	the kernel method is employed.
\end{abstract}


\maketitle

\section{Introduction}
Walter Kn\"odel introduced the following online bin-packing problem~\cite{walther}: There are bins of size 1, and
random items of size $\frac23$ (large items) and of size $\frac13$ (small items) appear and are put into the boxes.
A typical scenario is that  a number $j$ of partially filled boxes exist, and the number $j$ becomes $j+1$ resp.\ 
$j-1$, depending on whether a the new item is of large resp.\ small type. ``At random'' means that both types appear
with the same probability $\frac12$.

In my collection of examples \cite{Prodinger-kernel}, I showed how to deal with the Kn\"odel problem using the
kernel method. I was, however, not the only author who was intrigued by such questions; a notable paper is
by Michael Drmota \cite{michi}, which is of a more probabilistic type, whereas I tried to emphasize the combinatorial
point of view.

The present paper has a certain `ternary' flavour: the next section deals with the instance of large items appearing
in double-packs. The handler breaks off the double-packs, and then treats the items as Kn\"odel would have done.
Typically, the number of partially filled boxes increases by 2 or decreases by 1. In order to keep the system
balanced, we assume that the small items appear twice as often as the double-packs.

The last section deals with the dual problem, where the small items appear in double-packs and the large items as
single units.

The kernel method is used to obtain all the relevant enumerations. The recent paper \cite{Prodinger-PUMA} served
as an inspiration, but deals with a different issue. It must be said that, when \cite{Prodinger-kernel} was prepared,
such ternary questions would have been outside of my reach. Luckily, now, they are not.

We confine ourselves here just to enumerations, deriving explicit generating functions in one or two variables.
Questions of a more probabilistic nature are not treated.

\section{The first model}

The following items arrive at random: a double-pack of items, each of size $\frac23$, and an item of size $\frac13$.
We could equip the set-up with general probabilities $p$ and $q=1-p$, but we restrict ourselves to the `balanced' case
where the single items are twice as likely as the double-packs, so we set $p=\frac13$ and $q=\frac23$.

The following state diagram (we show only a finite part of it) describes the situation. There are states representing 
`$i$ boxes filled to $\frac23$'; a double-pack pushes the $i$ to $i+2$, and a single item reduces it to $i-1$. There is
an exceptional state, called $\beta$, standing for one box, filled to $\frac13$.
The red edges represent an arrival of a double-pack, and will be labelled by $pz$;
the black edges represent an arrival of a single item, and will be labelled by $qz$.

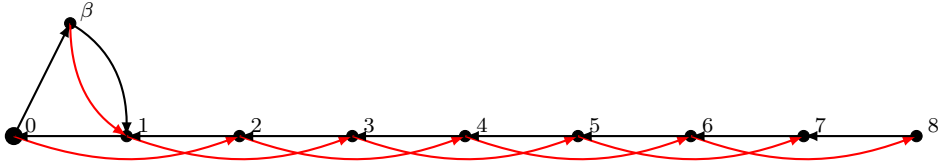
\begin{figure}[h]

\begin{center}
\begin{tikzpicture}[scale=1.5]

	\foreach \x in {0,1,2,3,4,5,6,7,8}
	{
		\draw (\x,0) circle (0.05cm);
	\fill (\x,0) circle (0.05cm);
	}

\draw (0.5,1) circle (0.05cm);	\fill (0.5,1) circle (0.05cm);
\fill (0,0) circle (0.08cm);

\foreach \x in {0,1,2,3,4,5,6}
{
\draw[thick, red, -latex] (\x,0) to [out=-20,in=200] (\x+2,0);	

}

\foreach \x in {0,1,2,3,4,5,6,7}
{
	\draw[thick, -latex] (\x+1,0) to  (\x,0);	
	\node at  (\x+0.15,0.1){\tiny$\x$};
}
\node at  (8+0.15,0.1){\tiny$8$};
\node at  (0.65,1.1){\tiny$\beta$};

\draw[thick, -latex] (0,0) to  (0.5,1);	
\draw[thick,red, -latex]   (0.5,1) [out=-90,in=145] to (1,0);
\draw[thick, -latex]   (0.5,1) [out=-30,in=90] to (1,0);	

\end{tikzpicture}
\end{center}
\end{figure}
From the state diagram, we set off an infinite set of generating functions in the variable $z$, where the coefficient of $z^n$ is the probability that $n$ random steps lead to state $i$, for $i\ge0$ or $i=\beta$. Mostly, we just write $f_i$ instead of $f_i(z)$. The following system of recursions can be read off immediately:
\begin{align*}
f_0&=1+qzf_1,\quad
f_\beta=qzf_0,\\
f_1&=zf_\beta+qzf_2=qz^2f_0+qzf_2,\\
f_i&=pzf_{i-2}+qzf_{i+1},\quad i\ge2.
\end{align*}
Our method to solve this system is the kernel method. For that, we introduce a bivariate generating function $F(u,z)$, but we mostly write just $F(u)$:
\begin{align*}
F(u)&=\sum_{i\ge0}u^if_i(z)\\
&=1+qzf_1+qz^2uf_0+qzuf_2+\sum_{i\ge2}u^i\Big[pzf_{i-2}+qzf_{i+1}\Big]\\
&=1+qzf_1+qz^2uf_0+pzu^2F(u)+\sum_{i\ge1}u^iqzf_{i+1}\\
&=1+qzf_1+qz^2uf_0+pzu^2F(u)+\frac{qz}{u}\Big(F(u)-f_0-uf_1\Big)\\
&=1+qz^2uf_0+pzu^2F(u)+\frac{qz}{u}\Big(F(u)-f_0\Big)\\
\end{align*}
Note that $f_0=F(0)$. It is beneficial to introduce the new variable $u=zU$; doing this, powers of $z$ that appear are  multiples of $3$. As can be seen, the numbers of steps leading to a state $i$ belong to just one residue class modulo 3. We compute
\begin{equation*}
F(u)=\frac{-3U-2z^3U^2f_0+2f_0}{z^3U^3-3U+2}=\frac{-3U-2xU^2f_0+2f_0}{xU^3-3U+2}.
\end{equation*}
As it is common using the kernel method, setting $U=0$ leads to a void equation. However, factorizing the denominator is the method of choice.
There is `bad' factor in the denominator, which must also appear in the numerator, which allows us to compute $f_0$ and consequently the whole
bivariate generating function. In order to deal with the ternary equation successfully, we further set
 $x=z^3=\frac{27}{4}t(1-t)^2$
and we find the 3 roots
\begin{align*}
U_1&=\frac2{3(1-t)},\quad
U_2=\frac{1}{\sigma},\quad
U_3=\frac{1}{\tau},
\end{align*}
with
\begin{align*}
\sigma&=\frac34(t-\sqrt{4t-3t^2}\,),\quad
\tau=\frac34(t+\sqrt{4t-3t^2}).
\end{align*}
Plugging $U=\frac2{3(1-t)}$ into the numerator (this is the bad factor, as explained a little bit later), leads to
\begin{equation*}
f_0=\frac{1}{(1-t)(1-3t)}
\end{equation*}
and furthermore to the simplified numerator
\begin{equation*}
\frac{-3U-2xU^2f_0+2f_0}{U-\frac2{3(1-t)}}=\frac1{1-3t}\Big(-3+\frac{27}{2}t(t-1)U\Big).
\end{equation*}
The variable $x$ is given in terms of $t$. The inverse relation is of interest. It can be obtained by
the Lagrange inversion formula or, as here, by contour integration:
\begin{align*}
[x^k]t&=\frac{1}{2\pi i}\oint\frac{dx}{x^{k+1}}t=
\frac{1}{2\pi i}\frac{27}{4}\Big(\frac{4}{27}\Big)^{k+1}\oint\frac{dt(1-t)(1-3t)}{t^{k+1}(1-t)^{2k+2}}t\\
&=\frac{1}{2\pi i}\Big(\frac{4}{27}\Big)^{k}\oint\frac{dt(1-3t)}{t^{k}(1-t)^{2k+1}}
=\Big(\frac{4}{27}\Big)^{k}[t^{k-1}]\frac{1-3t}{(1-t)^{2k+1}}\\
&=\Big(\frac{4}{27}\Big)^{k}\bigg[\binom{3k-1}{k-1}-3\binom{3k-2}{k-2}\bigg],
\end{align*}
which, after simplification, gives us
\begin{equation*}
	t=\sum_{k\ge1}\frac1k\binom{3k-2}{k-1}\frac{2^{2k}}{3^{3k}}.
\end{equation*}
A similar computation leads to
\begin{equation*}
	\frac1{1-t}=\sum_{k\ge0}\frac1{2k+1}\binom{3k}{k}\frac{2^{2k+1}}{3^{3k+1}}.
\end{equation*}
From this we infer that $U\sim \frac23$, or $u\sim\frac23z$, explaining why we are talking about the bad factor.
We continue the computation:
\begin{align*}
	F(u)&=\frac1{1-3t}\Big(-3+\frac{27}{2}t(t-1)U\Big)\frac1{x(U-\frac{1}{\sigma})(U-\frac{1}{\tau})}\\
&=\frac1{1-3t}\Big(-3+\frac{27}{2}t(t-1)U\Big)\frac{\frac94t(t-1)}{x(1-\sigma U)(1-\tau U)}\\
	&=\frac1{(1-3t)(1-t)}\Big(1-\frac{9}{2}t(t-1)U\Big)\frac{1}{(1-\sigma U)(1-\tau U)}.
\end{align*}
Partial fraction decomposition leads to (we use the abbreviation $W=\sqrt{4t-3t^2}\,$)
\begin{align*}
		\frac{1}{(1-\sigma U)(1-\tau U)}&=\frac12\Big(1-\frac tW\Big)\frac{1}{1-\sigma U}+\frac12\Big(1+\frac tW\Big)\frac{1}{1-\tau U}\\
		&=\frac12\Big[\frac{1}{1-\sigma U}+\frac{1}{1-\tau U}\Big]
		+\frac{t}{2W}\Big[\frac{1}{1-\tau U}-\frac{1}{1-\sigma U}\Big]\\
		&=\frac12\Big[\frac{1}{1-\sigma U}+\frac{1}{1-\tau U}\Big]
+\frac{3t}{4(\tau-\sigma)}\Big[\frac{1}{1-\tau U}-\frac{1}{1-\sigma U}\Big].
\end{align*}
For the further simplification we will resort to two identities
going by the name of Girard-Waring formula, see e. g. \cite{Gould}:
\begin{equation*}
	X^m+Y^m=\sum_{0\le k\le m/2}(-1)^k\binom{m-k}{k}\frac{m}{m-k}(XY)^{k}(X+Y)^{m-2k};
\end{equation*}
\begin{equation*}
	\frac{X^m-Y^m}{X-Y}=\sum_{0\le k\le (m-1)/2}(-1)^k \binom{m-1-k}{k}(XY)^{k}(X+Y)^{m-1-2k}.
\end{equation*}
Of course, we will apply them with $X=\tau$ and $Y=\sigma$. Then
\begin{align*}
	[U^m]\frac12\Big[\frac{1}{1-\sigma U}+\frac{1}{1-\tau U}\Big]
	&=\frac12\sum_{0\le k\le m/2}(-1)^k\binom{m-k}{k}\frac{m}{m-k}\Big(\frac94t(t-1)\Big)^{k}\Big(\frac32t\Big)^{m-2k}\\
	&=\frac12\Big(\frac32\Big)^{m}\sum_{0\le k\le m/2}(-1)^k\binom{m-k}{k}\frac{m}{m-k} t^k(t-1)^{k}t^{m-2k}\\
		&=\frac12\Big(\frac32\Big)^{m}\sum_{0\le k\le m/2}(-1)^k\binom{m-k}{k}\frac{m}{m-k} (t-1)^{k}t^{m-k}
\end{align*}
and
\begin{align*}
[U^m]&\frac{3t}{4(\tau-\sigma)}\Big[\frac{1}{1-\tau U}-\frac{1}{1-\sigma U}\Big]\\&
=\frac{3t}{4}\sum_{0\le k\le (m-1)/2}(-1)^k \binom{m-1-k}{k}\Big(\frac94t(t-1)\Big)^{k}\Big(\frac32t\Big)^{m-1-2k}\\
&=\Big(\frac32\Big)^{m-1}\frac{3t}{4}\sum_{0\le k\le (m-1)/2}(-1)^k \binom{m-1-k}{k}(t-1)^{k}t^{m-1-k}\\
&=\frac{1}{2}\Big(\frac32\Big)^{m}\sum_{0\le k\le (m-1)/2}(-1)^k \binom{m-1-k}{k}(t-1)^{k}t^{m-k}.
\end{align*}
Combining the two leads to a pleasant simplification:
\begin{align*}
[U^m]&\frac12\Big[\frac{1}{1-\sigma U}+\frac{1}{1-\tau U}\Big]+
[U^m]\frac{3t}{4(\tau-\sigma)}\Big[\frac{1}{1-\tau U}-\frac{1}{1-\sigma U}\Big]\\
&=\Big(\frac32\Big)^{m}\sum_{0\le k\le m/2}(-1)^k \binom{m-k}{k}(t-1)^{k}t^{m-k},
\end{align*}
or simpler
\begin{align*}
[U^m]\frac{1}{(1-\sigma U)(1-\tau U)}=\Big(\frac32\Big)^{m}\sum_{0\le k\le m/2}(-1)^k \binom{m-k}{k}(t-1)^{k}t^{m-k}.
\end{align*}
We need a second similar term:
\begin{align*}
	[U^m]&\Big(-\frac{9}{2}t(t-1)U\Big)\frac{1}{(1-\sigma U)(1-\tau U)}\\
	&=-\frac{9}{2}t(t-1)[U^{m-1}]\frac{1}{(1-\sigma U)(1-\tau U)}\\
	&=-\frac{9}{2}t(t-1)\Big(\frac32\Big)^{m-1}\sum_{0\le k\le m/2}(-1)^k \binom{m-1-k}{k}(t-1)^{k}t^{m-1-k}\\
		&=-3\Big(\frac32\Big)^{m}\sum_{0\le k\le m/2}(-1)^k \binom{m-1-k}{k}(t-1)^{k+1}t^{m-k}.
\end{align*}
Putting everything together we found
\begin{align*}
	F(u)
	&=\frac1{(1-3t)(1-t)}\sum_{m\ge0}\frac{u^m}{z^m}\Big(\frac32\Big)^{m}\sum_{0\le k\le m/2}(-1)^k \binom{m-k}{k}(t-1)^{k}t^{m-k}\\
	&-\frac3{(1-3t)(1-t)}\sum_{m\ge0}\frac{u^m}{z^m}\Big(\frac32\Big)^{m}\sum_{0\le k\le m/2}(-1)^k \binom{m-1-k}{k}(t-1)^{k+1}t^{m-k}.
	\end{align*}
Reading off coefficients of powers of $z$ as well is now done with Cauchy's integral formula; the contours are always small circles (or equivalent) around the origin.
The starting point is
\begin{align*}
	[z^nu^j]F(u)
&=\Big(\frac32\Big)^{j}[z^{n+j}]\frac1{(1-3t)(1-t)}\sum_{0\le k\le j/2}(-1)^k \binom{j-k}{k}(t-1)^{k}t^{j-k}\\*
	&-\Big(\frac32\Big)^{j}[z^{n+j}]\frac3{(1-3t)(1-t)}\sum_{0\le k\le j/2}(-1)^k \binom{j-1-k}{k}(t-1)^{k+1}t^{j-k}
\end{align*}
and we will treat the two sums separately. 
There is only a contribution if $n+j\equiv 0\bmod 3$. (This is also clear from the combinatorial context.) Assume this and set $N:=\frac{n+j}{3}$.

Step 1: 
\begin{align*}
[x^{N}]&\frac1{(1-3t)(1-t)}\sum_{0\le k\le j/2}(-1)^k \binom{j-k}{k}(t-1)^{k}t^{j-k}\\
&=\frac{1}{2\pi i}\oint\frac{dx}{x^{N+1}}\frac1{(1-3t)(1-t)}\sum_{0\le k\le j/2}(-1)^k \binom{j-k}{k}(t-1)^{k}t^{j-k}\\
&=\frac{1}{2\pi i}\oint\frac{27}{4}\frac{dt}{\big(\frac{27}{4}t(1-t)^2\big)^{N+1}}\sum_{0\le k\le j/2}(-1)^k \binom{j-k}{k}(t-1)^{k}t^{j-k}\\
&=\frac{1}{2\pi i}\oint\Big(\frac{4}{27}\Big)^{N}dt\sum_{0\le k\le j/2}(-1)^k \binom{j-k}{k}(t-1)^{k-2N-2}t^{j-k-N-1}\\
&=\Big(\frac{4}{27}\Big)^{N}[t^{N-j+k}]\sum_{0\le k\le j/2}(-1)^k \binom{j-k}{k}(t-1)^{k-2N-2}\\
&=\Big(\frac{4}{27}\Big)^{N}(-1)^{N-j}\sum_{0\le k\le j/2} \binom{j-k}{k}\binom{k-2N-2}{N-j+k}.
\end{align*}
Step 2: 
\begin{align*}
	[x^{N}]&\frac3{(1-3t)(1-t)}\sum_{0\le k\le j/2}(-1)^k \binom{j-1-k}{k}(t-1)^{k+1}t^{j-k}\\
	&=3\Big(\frac{4}{27}\Big)^{N}[t^{N-j+k}]\sum_{0\le k\le j/2}(-1)^k \binom{j-1-k}{k}(t-1)^{k-2N-1}\\
	&=3\Big(\frac{4}{27}\Big)^{N}(-1)^{N-j}\sum_{0\le k\le j/2} \binom{j-1-k}{k}\binom{k-2N-1}{N-j+k}.
\end{align*}
We put all the results of this section together in a theorem.

\begin{theorem}
	The generating function $F(u)=F(u,z)$ has the following explicit form:
\begin{equation*}
	F(u)=\frac1{(1-3t)(1-t)}\Big(1-\frac{9}{2}t(t-1)U\Big)\frac{1}{(1-\sigma U)(1-\tau U)}
\end{equation*}
Here, $u=zU$, $z^3=x=\frac{27}{4}t(1-t)^2$, and
\begin{align*}
	\sigma&=\frac34(t-\sqrt{4t-3t^2}\,),\quad
	\tau=\frac34(t+\sqrt{4t-3t^2}\,).
\end{align*}
Note that $(1-\sigma U)(1-\tau U)=1-\frac32tU+\frac94t(t-1)U^2$.
	Written in the new variable $U$, only powers of $z$ that are multiples of $3$ appear. Further, we get
	the representation sorted by powers of $u$:
\begin{align*}
	F(u,z)
	&=\frac1{(1-3t)(1-t)}\sum_{m\ge0}\frac{u^m}{z^m}\Big(\frac32\Big)^{m}\sum_{0\le k\le m/2}(-1)^k \binom{m-k}{k}(t-1)^{k}t^{m-k}\\
	&-\frac3{(1-3t)(1-t)}\sum_{m\ge0}\frac{u^m}{z^m}\Big(\frac32\Big)^{m}\sum_{0\le k\le (m-1)/2}(-1)^k \binom{m-1-k}{k}(t-1)^{k+1}t^{m-k}.
\end{align*}
Reading off coefficients of $z^Nu^j$, where $N=\frac{n+j}{3}$ leads to
\begin{align*}
[z^Nu^j]F(u,z)&=\Big(\frac{4}{27}\Big)^{N}(-1)^{N-j}\sum_{0\le k\le j/2} \binom{j-k}{k}\binom{k-2N-2}{N-j+k}\\
&-3\Big(\frac{4}{27}\Big)^{N}(-1)^{N-j}\sum_{0\le k\le (j-1)/2} \binom{j-1-k}{k}\binom{k-2N-1}{N-j+k}.
\end{align*}
	For the special state $\beta$, the following series representation holds:
	\begin{align*}
		f_\beta(z)=\sum_{n\ge0}\frac{2^{2n+1}}{{3}^{3n+1}}\binom{3n+1}{n}z^{3n+1}.
	\end{align*}\qed
\end{theorem}

The computation for the special state was not shown yet:
\begin{align*}
	[z^{3n+1}]&f_\beta=\frac23[x^n]\frac1{(1-t)(1-3t)}
	=\frac23\frac1{2\pi i}\oint\frac{dx}{x^{n+1}}\frac1{(1-t)(1-3t)}\\
	&=\frac23\frac{27}{4}\frac1{2\pi i}\oint\frac{dt}{(\tfrac{27}4)^{n+1}t^{n+1}(1-t)^{2n+2}}
	=\frac23 \Big(\frac{4}{27}\Big)^n\frac1{2\pi i}\oint\frac{dt}{t^{n+1}(1-t)^{2n+2}}\\
	&=\frac23 \Big(\frac{4}{27}\Big)^n[t^n]\frac{1}{(1-t)^{2n+2}}=\frac{2^{2n+1}}{{3}^{3n+1}}\binom{3n+1}{n}.
\end{align*}


\section{The dual model}

Now, the red edges mean the arrival of the large objects (size $\frac23$) and the black edges mean a double-pack of the small edges
(size $\frac13$ each). To keep the system balanced, the large objects should arrive twice as often as the double-packs of small edges. 
Again, the generating function $g_i$ refers to paths of length $n$ leading eventually into state $i$. After $n$ steps, only a state $i$ can be
reached with $n\equiv i \bmod 3$. The state diagram and the recursions are immediate:
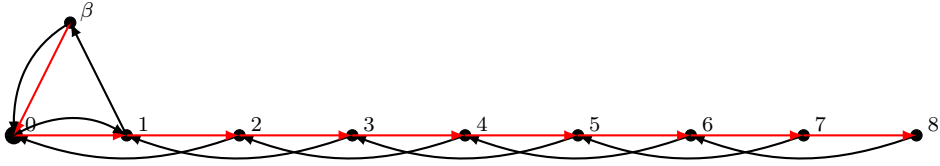
\begin{figure}[h]

	\begin{center}
		\begin{tikzpicture}[scale=1.5]

			\foreach \x in {0,1,2,3,4,5,6,7,8}
			{
				\draw (\x,0) circle (0.05cm);
				\fill (\x,0) circle (0.05cm);
			}
			
			\draw (0.5,1) circle (0.05cm);	\fill (0.5,1) circle (0.05cm);
			\fill (0,0) circle (0.08cm);

			\foreach \x in {0,1,2,3,4,5,6}
			{
				\draw[thick,  latex-] (\x,0) to [out=-20,in=200] (\x+2,0);	
				
			}
			
			\foreach \x in {0,1,2,3,4,5,6,7}
			{
				\draw[thick, red,latex-] (\x+1,0) to  (\x,0);	
				\node at  (\x+0.15,0.1){\tiny$\x$};
			}
			\node at  (8+0.15,0.1){\tiny$8$};
			\node at  (0.65,1.1){\tiny$\beta$};
			
			
			\draw[thick, -latex] (0,0)[out=30,in=150] to  (1,0);
			\draw[thick, latex-]   (0.5,1)  to (1,0);
			\draw[thick, red,-latex]   (0.5,1)   to (0,0);
			\draw[thick, -latex]   (0.5,1) [out=210,in=90]  to (0,0);	
		\end{tikzpicture}
	\end{center}
\end{figure}

We work only with $p=\frac23$, $q=\frac13$. Directly from the state diagram,
\begin{align*}
g_0&=1+zg_\beta+qzg_2=1+qz^2g_1+qzg_2,\\
g_\beta&=qzg_1, \quad
g_1=zg_0+qzg_3,\\
g_i&=pzg_{i-1}+qzg_{i+2},\quad i\ge2.
\end{align*}
Summing the recursions,
\begin{align*}
G(u)&=g_0+ug_1+\sum_{i\ge2}u^i\Big(pzg_{i-1}+qzg_{i+2}\Big)\\
&=g_0+uzg_0+qzug_3+pzu\sum_{i\ge1}u^ig_i+\frac{qz}{u^2}\sum_{i\ge4}u^ig_{i}\\
&=g_0+uzg_0+pzuG(u)-pzug_0+\frac{qz}{u^2}\sum_{i\ge3}u^ig_{i}\\
&=g_0+uzg_0+pzuG(u)-pzug_0+\frac{qz}{u^2}(G(u)-g_0-ug_1-u^2g_2)\\
&=g_0+quzg_0+pzuG(u)+\frac{qz}{u^2}G(u)-\frac{qz}{u^2}g_0-\frac{qz}{u}g_1-qzg_2\\
&=g_0+quzg_0+pzuG(u)+\frac{qz}{u^2}G(u)-\frac{qz}{u^2}g_0-\frac{qz}{u}g_1+1+qz^2g_1-g_0.
\end{align*}
Solving, we find with $V=uz$:
\begin{equation*}
G(u)=\frac{-{V}^{3}g_0-3{V}^{2}-g_1{V}^{2}{z}^{2}+{z}^{3}g_0+g_1V{z}^{2}}{2V^3-3V^2+x}.
\end{equation*}
Now we factorize the denominator:
\begin{equation*}
2(V-\tfrac32(1-t))(V-\sigma)(V-\tau)=2V^3-3V^2+x.
\end{equation*}
This time, both, $(V-\sigma)$ and $(V-\tau)$ are bad factors. Plugging into the numerator, we find two equations, and the solutions:
\begin{equation*}
g_0=\frac{4}{(1-3t)(4-3t)},\quad g_1=\frac{27t(1-t)}{z^2(1-3t)(4-3t)}.
\end{equation*}
Dividing out the bad factors, we find
\begin{equation*}
\frac{-{V}^{3}g_0-3{V}^{2}-g_1{V}^{2}{z}^{2}+{z}^{3}g_0+g_1V{z}^{2}}{(V-\sigma)(V-\tau)}=
\frac{12(t-1)-4V}{(1-3t)(4-3t)}.
\end{equation*}
Altogether:
\begin{align*}
G(u)&=		\frac{6(t-1)-2V}{(1-3t)(4-3t)}\frac{1}{V-\frac32(1-t)}\\
		&=		\frac{6(1-t)+2V}{(1-3t)(4-3t)\frac32(1-t)}\frac{1}{1-\frac2{3(1-t)}V}\\
		&=		2\frac{2+\frac2{3(1-t)}V}{(1-3t)(4-3t)}\frac{1}{1-\frac2{3(1-t)}V}.
\end{align*}
Furthermore
\begin{align*}
[V^j]G(u)&=\frac2{(1-3t)(4-3t)}\bigg[2\Big(\frac23\frac{1}{1-t}\Big)^{j}+\Big(\frac23\frac{1}{1-t}\Big)^{j}\bigg]\\
&=\frac6{(1-3t)(4-3t)}\Big(\frac23\frac{1}{1-t}\Big)^{j},\quad j\ge1,
\end{align*}
and
\begin{align*}
	[u^j]G(u)=z^j\frac6{(1-3t)(4-3t)}\Big(\frac23\frac{1}{1-t}\Big)^{j}.
\end{align*}
Now let us consider $j+3N$ steps to reach state $j$, and then
\begin{align*}
	[z^{j+3N}u^j]G(u)&=[x^N]\frac6{(1-3t)(4-3t)}\Big(\frac23\frac{1}{1-t}\Big)^{j}\\
&=\frac1{2\pi i}\oint\frac{dx}{x^{N+1}}	\frac6{(1-3t)(4-3t)}\Big(\frac23\frac{1}{1-t}\Big)^{j}\\
&=\frac{27}{4}\Big(\frac4{27}\Big)^{N+1}\frac1{2\pi i}\oint\frac{dt(1-t)(1-3t)}{t^{N+1}(1-t)^{2N+2}}	\frac6{(1-3t)(4-3t)}\Big(\frac23\frac{1}{1-t}\Big)^{j}\\
&=\Big(\frac4{27}\Big)^{N}\Big(\frac23\Big)^{j}\frac1{2\pi i}\oint\frac{dt}{t^{N+1}(1-t)^{2N+j+1}}	\frac6{(4-3t)}\\
&=\Big(\frac4{27}\Big)^{N}\Big(\frac23\Big)^{j-1}[t^{N}]\frac{1}{(1-t)^{2N+j+1}}	\frac1{(1-\frac34t)}\\
&=\Big(\frac4{27}\Big)^{N}\Big(\frac23\Big)^{j-1}\sum_{i=0}^N\Big(\frac34\Big)^{N-i}\binom{2N+j+i}{i}	\\
&=\sum_{i=0}^N\frac{2^{2i+j-1}}{3^{2N+i+j-1}}\binom{2N+j+i}{i}.
\end{align*}
The coefficients of $g_0$ are different:
\begin{equation*}
[z^{3N}]g_0=\sum_{i=0}^N\frac{2^{2i}}{3^{2N+i}}\binom{2N+i}{i}.
\end{equation*}
Furthermore,
\begin{equation*}
	[z^{3N+1}]g_\beta=\frac13[z^{3N}]g_1=\sum_{i=0}^N\frac{2^{2i}}{3^{2N+i+1}}\binom{2N+1+i}{i}.
\end{equation*}

Here are the main results of this section:

\begin{theorem}
	The generating function $G(u)=G(u,z)$ has the following explicit form:
	\begin{equation*}
		G(u)=			2\frac{2+\frac2{3(1-t)}V}{(1-3t)(4-3t)}\frac{1}{1-\frac2{3(1-t)}V}.
	\end{equation*}
	Here, $u=\frac Vz$, $z^3=x=\frac{27}{4}t(1-t)^2$.
	Written in the new variable $V$, only powers of $z$ that are multiples of $3$ appear. Further, we get
	the representation sorted by powers of $u$:
\begin{align*}
	[V^j]G(u)
	&=\frac6{(1-3t)(4-3t)}\Big(\frac23\frac{1}{1-t}\Big)^{j},\quad j\ge1,
\end{align*}
and
\begin{align*}
	[u^j]G(u)=z^j\frac6{(1-3t)(4-3t)}\Big(\frac23\frac{1}{1-t}\Big)^{j}.
\end{align*}
	Reading off coefficients of $z^{j+3N}u^j$ leads to
	\begin{align*}
		[z^{j+3N}u^j]G(u,z)=\sum_{i=0}^N\frac{2^{2i+j-1}}{3^{2N+i+j-1}}\binom{2N+j+i}{i}.
	\end{align*}
	For the special cases, the following series representation holds:
	\begin{align*}
		[z^{3N}]g_0&=\sum_{i=0}^N\frac{2^{2i}}{3^{2N+i}}\binom{2N+i}{i},\\
		[z^{3N+1}]g_\beta&=\frac13[z^{3N}]g_1=\sum_{i=0}^N\frac{2^{2i}}{3^{2N+i+1}}\binom{2N+1+i}{i}.
	\end{align*}\qed
\end{theorem}

\newpage

\bibliographystyle{plain}


\end{document}